\documentclass[11pt]{article}
\usepackage[T1]{fontenc}
\usepackage{amsmath,amsfonts,amsthm,amssymb,mathtools}
\usepackage{mathrsfs,graphicx,color,latexsym,calc}
\usepackage[colorlinks,bookmarksopen,bookmarksnumbered,citecolor=blue,linkcolor=red,urlcolor=blue]{hyperref}
\usepackage{enumitem}
\usepackage{authblk}
\usepackage{booktabs,array}

\allowdisplaybreaks[2]
\newtheorem{theorem}{Theorem}[section]
\newtheorem{lemma}[theorem]{Lemma}

\newtheorem{conjecture}{Conjecture}[section]
\newtheorem*{conjecture*}{Conjecture}
\theoremstyle{definition}

\theoremstyle{remark}
\newtheorem{remark}[theorem]{Remark}
\numberwithin{equation}{section}

\newcommand{\one}{\mathbf{1}}
\newcommand{\EE}{\mathbb{E}}
\newcommand{\PP}{\mathbb{P}}
\newcommand{\Zge}{\mathbb{Z}_{\ge0}}

\hypersetup{
  pdftitle={A near-linear upper bound for Burr’s conjecture},
  pdfauthor={Liangdong Fan, Junying Lu, Yaojun Chen},
  pdfsubject={Absorbing sets of small chromatic number and oriented tree embeddings},
  pdfkeywords={oriented tree, chromatic number, absorbing set, Burr conjecture, Lovasz local lemma}
}

\begin{document}

\title{A near-linear upper bound for Burr’s conjecture}
% Replace the empty author field before submission.  For example:
% \author[1]{First Author}
% \affil[1]{Department, University, City, Country}
\author{Liangdong FAN$^1$, Junying LU$^2$, Yaojun CHEN$^{1,}$\footnote{Corresponding author. Email: yaojunc@nju.edu.cn}\\
 {\small{$^1$School of Mathematics, Nanjing University, Nanjing 210093, P.R. CHINA}}\\
 {$^2$\small{School of Mathematics and Statistics, Nanjing University of Science and Technology,\\
Nanjing 210094, P.R. CHINA}}\\
 }
\date{}
\maketitle

\begin{abstract}
Let $f(k)$ denote the smallest integer such that every oriented graph $D$ with chromatic number at least $f(k)$ contains every
oriented tree on $k$ vertices.
Burr (1980) showed that $f(k)\le (k-1)^2$ and conjectured that $f(k)=2k-2$. Bessy, Gon\c{c}alves and
Reinald (2025) proved that
$f(k)=O(k^{3/2})$. In this paper, by using an absorbing set method, we show that 
$f(k)\le \lfloor 31\log (k!)\rfloor=O(k\log k)$.

\vskip 2mm
\noindent\textit{\bf AMS classification:} 05C20; 05C15; 05C05; 05D40\\[1mm]
\noindent\textit{\bf Keywords:} Burr's conjecture; Oriented tree; Chromatic number; Absorbing set
\end{abstract}

\baselineskip=0.202in

\section{Introduction}
\label{sec:intro}
All digraphs considered in this paper are finite without loops,
parallel arcs or digons, referred to as oriented graphs. Let $D$ be an oriented graph with vertex set $V(D)$. The underlying graph of $D$ is an undirected graph obtained by neglecting its orientation. Let $\chi(D)$ denote the chromatic number of its underlying graph.
For $U\subseteq V(D)$, write $D[U]$ for the subgraph induced by $U$  and $D-U=D[V(D)\setminus U]$. The out- and in-neighbourhoods of $v\in V(D)$
are $N_D^+(v)$ and $N_D^-(v)$, respectively. Set $N_D^\pm[v]=N_D^\pm(v)\cup\{v\}$,
$d_D^\pm(v)=|N_D^\pm(v)|$, and
$\Delta^\pm(D)=\max_{v\in V(D)}d_D^\pm(v)$. If $V(D)=\varnothing$, we set $\chi(D)=\Delta^\pm(D)=0$.
%We put $\chi(\varnothing)=\Delta^\pm(\varnothing)=0$.  
We use $\log$ and $\ln$ to denote the logarithmic function with bases 2 and $e$, respectively. 

It is well known that every undirected graph of chromatic number at least $k$ contains every tree on $k$
vertices: a vertex-critical subgraph of chromatic number at least $k$ has
minimum degree at least $k-1$, and the tree can then be embedded greedily.
For oriented trees, this argument breaks down, since it gives no control over
the directions of the unused edges. Let $f(k)$ denote the smallest integer such that every oriented graph $D$ with $\chi(D)\ge f(k)$ contains every
oriented tree on $k$ vertices. Burr \cite{Burr1980}
proved that $f(k)\le(k-1)^2$ and posed the following.
\begin{conjecture}[Burr \cite{Burr1980}]\label{c1}
$f(k)=2k-2$.
\end{conjecture}
%For a digraph $D$, let $\chi(D)$ denote the chromatic number of its underlying undirected graph. 
%An oriented tree is an orientation of a tree, and all copies
%considered here are non-induced. For $k\ge2$, let $b(k)$ be the least integer
%$c$ such that every oriented graph $D$ with $\chi(D)\ge c$ contains every
%oriented tree on $k$ vertices. It is immaterial whether the host is required
%to be oriented: deleting one arc from each digon still preserves the underlying
%graph and hence its chromatic number.
The lower bound $f(k)\ge 2k-2$ follows
from a regular tournament of order $2k-3$, in which every vertex has
out-degree $k-2$ and hence the out-star on $k$ vertices cannot occur. 
%Burr proved the upper bound $b(k)\le(k-1)^2$. 
Many years later, Addario-Berry, Havet, Linhares Sales,
Reed and Thomass\'e~\cite{AddarioBerryEtAl2013} improved this to
$f(k)\le k(k-1)/2+1$. The first subquadratic
estimate is due to 
Bessy, Gon\c{c}alves and
Reinald~\cite{BessyGoncalvesReinald2025} who showed the following.
 
\begin{theorem}[Bessy, Gon\c{c}alves and
Reinald~\cite{BessyGoncalvesReinald2025}] 
$f(k)\le8\sqrt{2/15}\,k^{3/2}+O(k)$.
\end{theorem}

Linear bounds are known under additional assumptions. The first one is the famous 
Gallai--Hasse--Roy--Vitaver theorem which gives the exact threshold $k$ for the
directed path on $k$ vertices~\cite{Gallai1968,Hasse1965,Roy1967,Vitaver1962}.
Addario-Berry et al.~\cite{AddarioBerryEtAl2013} proved linear bounds for
antidirected trees and showed that chromatic number $k$ suffices when the
host is acyclic. Most recently, Santos, Stein and Williams
\cite{SantosSteinWilliams} and, independently, Riordan and Scott
\cite{RiordanScott} settled Burr's conjecture for antidirected trees. A different, host-dependent formulation was studied by
Naia~\cite{Naia2022}, who proved that every orientation of a graph $G$
contains every oriented tree of order
$\lfloor\chi(G)/\log|V(G)|\rfloor$.

The corresponding problem for tournaments is the following well-known conjecture.

\begin{conjecture}[Sumner~\cite{ReidWormald1983}]
Every tournament on $2k-2$ vertices contains every oriented tree on $k$ vertices.
\end{conjecture}

%K\"uhn, Mycroft and Osthus~\cite{KuhnMycroftOsthus2011} proved the conjecture for all sufficiently large $k$.
 K\"uhn,
Mycroft and Osthus~\cite{KuhnMycroftOsthus2011} proved the conjecture for all
sufficiently large $k$. For arbitrary $k$, Dross and
Havet~\cite{DrossHavet2021} showed that
$\lceil21k/8-47/16\rceil$ vertices suffice. Recently, Ai and
Lian~\cite{AiLian2026} improved this to $\lceil(18k-23)/7\rceil$. Benford and
Montgomery~\cite{BenfordMontgomery2022} obtained the sharper bound $k+C\ell$
for trees with $\ell$ leaves, where $C$ is an absolute constant. These
results rely on the completeness of the host and do not directly yield
bounds for arbitrary digraphs in terms of chromatic number.
\vskip2mm
In this paper, we obtain a nearly linear bound without restrictions on the host or on the
orientation of the target tree. The main result is as follows.

\begin{theorem}
\label{thm:main}
Let $k\ge2$ be an integer and let $D$ be an oriented graph. If
$\chi(D)\ge \lfloor 31\log (k!)\rfloor$, then $D$ contains every oriented tree on $k$
vertices.
\end{theorem}

Clearly, Theorem \ref{thm:main} tells us $f(k)=O(k\log k)$, which improves the upper bound $O(k^{3/2})$ obtained by Bessy, Gon\c{c}alves and
Reinald~\cite{BessyGoncalvesReinald2025}. This, together with the regular tournament construction, narrows the
gap between the known upper and lower bounds to a factor $O(\log k)$.

A set $S\subseteq V(D)$ is called
\emph{absorbing} if every vertex of $V(D)\setminus S$ has an out-neighbour in
$S$. The proof of Theorem \ref{thm:main} rests on the following result concerning absorbing sets, which is of interesting in its own right. 

\begin{theorem}
\label{thm:absorb}
Let $d\ge0$ be an integer and let $D$ be an oriented graph with
$\Delta^+(D)\le d$. Then $D$ has an absorbing set $S$ such that
$\chi(D[S])\le31\log(d+2)$.
\end{theorem}

The point is that the bound is independent of both $\Delta^-(D)$ and
$|V(D)|$. Reserving $S$ therefore costs at most $31\log(d+2)$ colours, since
$\chi(D-S)\ge\chi(D)-\chi(D[S])$.

\begin{remark}
\label{rem:absorb-optimal}
The logarithmic dependence on $d$ in Theorem~\ref{thm:absorb} is best
possible up to a constant factor. Indeed, a fixed $s$-vertex set in a random
tournament $R_n$ is absorbing with probability $(1-2^{-s})^{n-s}$. Hence,
for every fixed $0<\varepsilon<1$, a union bound shows that, asymptotically
almost surely, $R_n$ has no absorbing set of order at most
$(1-\varepsilon)\log n$. Since every induced subgraph of $R_n$ has complete
underlying graph and $\Delta^+(R_n)\le n-1$, taking $d=n-1$ shows that the
logarithmic order in Theorem~\ref{thm:absorb} cannot be improved.
\end{remark}

\section{Proof of Theorem~\ref{thm:main}}
\label{sec:embedding}
We first deduce Theorem~\ref{thm:main} from
Theorem~\ref{thm:absorb}, whose proof will be given in Section~3.
\begin{proof}
We proceed by induction on \(k\), the order of the target tree. The case \(k=2\) is trivial. Let \(k\ge 3\), and assume the result holds for \(k-1\).

Let \(T\) be an oriented tree of order \(k\), and let \(D\) be an oriented graph with \(\chi(D)\ge \lfloor 31\log (k!)\rfloor\). Choose a leaf \(z\) of \(T\), let \(p\) be its neighbour, and put \(T'=T-z\). Reversing \(D\) and \(T\) if necessary, we may assume that \(p\to z\).

Let \(L=\{v\in V(D):d_D^+(v)\le k-2\}\). Since \(\Delta^+(D[L])\le k-2\), by Theorem~\ref{thm:absorb} we can obtain an absorbing set \(S\subseteq L\) in \(D[L]\) with \(\chi(D[S])\le31\log k\). Since \(\chi(D[S])\) is an integer, \(\chi(D[S])\le\lfloor31\log k\rfloor\). Hence
\[
\chi(D-S)\ge\chi(D)-\chi(D[S])
\ge \lfloor31\log(k!)\rfloor-\lfloor31\log k\rfloor
\ge \lfloor31\log((k-1)!)\rfloor,
\]
where the last inequality follows from \(\lfloor a+b\rfloor\ge\lfloor a\rfloor+\lfloor b\rfloor\).
By induction, \(T'\) has an embedding \(\varphi\) in \(D-S\). Let \(x=\varphi(p)\).

If \(x\in L\), then \(x\in L\setminus S\), and since \(S\) is absorbing in \(D[L]\), there is \(y\in S\) with \(x\to y\). If \(x\notin L\), then \(d_D^+(x)\ge k-1\), while at most \(k-2\) out-neighbours of \(x\) lie in \(\varphi(V(T'))\). Thus in either case there is \(y\notin\varphi(V(T'))\) with \(x\to y\).

Setting \(\varphi(z)=y\) extends \(\varphi\) to an embedding of \(T\) in \(D\), completing the induction and hence the proof.
\end{proof}

\section{Proof of Theorem~\ref{thm:absorb}}

It remains to prove Theorem~\ref{thm:absorb}. We begin by recalling three standard tools that will be used in the proof. The first is Brouwer's fixed-point theorem~\cite{Brouwer1911}, the second is Markov's inequality~\cite{Durrett2019}, and the third is the symmetric Lov\'asz local lemma \cite{ErdosLovasz1975}.

\begin{theorem}[Brouwer~\cite{Brouwer1911}]
\label{thm:brouwer}
Let \(K\) be a nonempty compact convex subset of \(\mathbb{R}^n\), and let
\(f:K\to K\) be continuous. Then there exists \(\mathbf{x}\in K\) such that
\[
f(\mathbf{x})=\mathbf{x}.
\]
\end{theorem}

\begin{lemma}[Markov~\cite{Durrett2019}]
\label{lem:markov}
Let \(X\) be a nonnegative random variable. Then, for every \(a>0\),
\[
\PP(X\ge a)\le \frac{\EE[X]}{a}.
\]
\end{lemma}
%We use the following form of the symmetric Lov\'asz local lemma.
\begin{lemma}[Erd\H{o}s and Lov\'asz~\cite{ErdosLovasz1975}]
\label{lem:lll}
Let \(\mathcal X\) be a finite family of mutually independent random
variables, and let \(A_1,\ldots,A_m\) be events. Suppose that, for each
\(i\in[m]\), the event \(A_i\) is determined by some
\(\mathcal X_i\subseteq\mathcal X\), and that
\[
  \PP(A_i)\le p,
  \qquad
  \bigl|\{j\ne i:\mathcal X_i\cap\mathcal X_j\ne\varnothing\}\bigr|\le q
\]
for every \(i\in[m]\). If \(ep(q+1)\le1\), then 
\(\PP(\bigcap_{i=1}^m \overline{A_i})>0\).
\end{lemma}

Our first step towards proving Theorem \ref{thm:absorb} is to establish a weak version under the assumption that 
both $\Delta^+(D)$ and $\Delta^-(D)$ are bounded. For any oriented graph $D$, Theorem \ref{thm:brouwer} helps us find
weights $x_v\in[0,1]$, one assigned to each vertex $v\in V(D)$, such that
$\sum_{u\in N_D^+[v]}x_u\geq 1$, while
$\sum_{u\in N_D^+(v)}x_u\leq 1$ whenever $x_v>0$. Based on this weight function, by the probabilistic method, we can show that $D$  has an absorbing
set inducing a digraph of logarithmically bounded out-degree, and hence of
logarithmically bounded chromatic number.

%\paragraph{Proof strategy for the absorbing-set theorem.}
%We first prove the absorbing-set theorem when both $\Delta^+(D)$ and
%$\Delta^-(D)$ are bounded. A fixed-point argument gives weights
%$x_v\in[0,1]$ such that $\sum_{u\in N_D^+[v]}x_u\ge1$ for every $v$, while
%$\sum_{u\in N_D^+(v)}x_u\le1$ whenever $x_v>0$. Independent rounding and
%the symmetric Lov\'asz local lemma (Lemma~\ref{lem:lll}) then give an absorbing
%set inducing a digraph of logarithmically bounded out-degree, and hence of
%logarithmically bounded chromatic number.

%Section~\ref{sec:bounded} proves the bounded-degree case,
%Section~\ref{sec:layering} removes the in-degree restriction, and
%Section~\ref{sec:embedding} completes the leaf-extension argument.

%\section{Absorbing sets with bounded in- and out-degrees}
%\label{sec:bounded}

%We begin by recalling several standard tools that will be used in the proof.
%Using these tools, we prove the following lemma for oriented graphs of bounded in- and out-degree.
\begin{lemma}
\label{lem:bounded}
Let \(d,h\in\Zge\), and let \(D\) be an oriented graph with
\(\Delta^+(D)\le d\) and \(\Delta^-(D)\le h\). Set
\(\tau=\ln\bigl(2e(d+1)(h+1)\bigr)\) and \(t=\lceil3\tau\rceil\). Then
\(D\) has an absorbing set \(J\) such that \(\chi(D[J])\le2t-1\).
\end{lemma}

\begin{proof}
Define \(\Phi:[0,1]^{V(D)}\to[0,1]^{V(D)}\) by
\(\Phi(\mathbf{x})_v=\max\{0,1-\sum_{u\in N_D^+(v)}x_u\}\) for every
\(v\in V(D)\), where \(\mathbf{x}=(x_v)_{v\in V(D)}\). Since $\Phi$ is a continuous self-map of the compact convex set
$[0,1]^{V(D)}$, Theorem \ref{thm:brouwer} guarantees the existence
of a fixed point $\mathbf{x}$ of $\Phi$. Thus, for every
$v\in V(D)$, $x_v=\max\left\{0,\,1-\sum_{u\in N_D^+(v)}x_u\right\}$.
Consequently,
\[
\sum_{u\in N_D^+[v]}x_u\ge1,
\qquad
x_v>0\ \Longrightarrow\
\sum_{u\in N_D^+(v)}x_u=1-x_v\le1.
\]

For each \(v\in V(D)\), let \(\xi_v\) be an independent \(\{0,1\}\)-valued
random variable such that
\(\PP(\xi_v=1)=p_v:=\min\{1,\tau x_v\}\), and set
\(J=\{v\in V(D):\xi_v=1\}\), so that \(\PP(v\in J)=p_v\) for every
\(v\in V(D)\).

For \(v\in V(D)\), let
\(E_v=\{N_D^+[v]\cap J=\varnothing\}\). If \(p_u=1\) for some
\(u\in N_D^+[v]\), then \(\PP(E_v)=0\). Otherwise
\(p_u=\tau x_u\) for every \(u\in N_D^+[v]\). Using the elementary
inequality \(1-y\le e^{-y}\) for \(y\ge0\), together with the independence
of the variables \(\xi_u\), we obtain
\[
\PP(E_v)
=\prod_{u\in N_D^+[v]}(1-p_u)
\le \exp\left(-\sum_{u\in N_D^+[v]}p_u\right)
=\exp\left(-\tau\sum_{u\in N_D^+[v]}x_u\right)
\le e^{-\tau}.
\]

For each \(v\) with \(x_v>0\), set
\(X_v=|N_D^+(v)\cap J|\) and \(G_v=\{X_v\ge t\}\). Since
\(G_v=\{3^{X_v}\ge3^t\}\), by Lemma~\ref{lem:markov}, we have 
\(\PP(G_v)\le3^{-t}\EE[3^{X_v}]\). By the independence of the variables
\(\xi_u\), we have
\(\EE[3^{X_v}]=\prod_{u\in N_D^+(v)}\EE[3^{\xi_u}]
=\prod_{u\in N_D^+(v)}(1+2p_u)\).
Using \(1+y\le e^y\) for \(y\ge0\), together with
\(p_u\le\tau x_u\) and
\(\sum_{u\in N_D^+(v)}x_u\le1\), we obtain
\[
\begin{aligned}
\PP(G_v)
&\le 3^{-t}\exp\left(2\sum_{u\in N_D^+(v)}p_u\right)
 \le 3^{-t}\exp\left(2\tau\sum_{u\in N_D^+(v)}x_u\right)\\
&\le \exp(2\tau-t\ln3)
 \le \exp\bigl(-(3\ln3-2)\tau\bigr)
 <e^{-\tau},
\end{aligned}
\]
where the last two inequalities follow from \(t\ge3\tau\) and
\(3\ln3-2>1\).

To apply Lemma~\ref{lem:lll}, consider all the events \(E_v\), together with
the events \(G_v\) for which \(x_v>0\). The event \(E_v\) depends only on
the variables \(\{\xi_u:u\in N_D^+[v]\}\), while \(G_v\) depends only on
\(\{\xi_u:u\in N_D^+(v)\}\). Hence each event depends on at most \(d+1\)
variables.

Conversely, for a fixed \(u\in V(D)\), the variable \(\xi_u\) can occur only
in \(E_u\) and in the events \(E_v\) and \(G_v\) with
\(v\in N_D^-(u)\). Since \(d_D^-(u)\le h\), the variable \(\xi_u\) occurs
in at most \(2h+1\) events. Therefore each event shares a variable with at
most \(2h(d+1)\) other events.

Since \(\PP(E_v)\le e^{-\tau}\) and \(\PP(G_v)<e^{-\tau}\), we may take
\(p=e^{-\tau}\) and \(q=2h(d+1)\) in Lemma~\ref{lem:lll}. Then
\[
ep(q+1)
=\frac{2h(d+1)+1}{2(d+1)(h+1)}
<1.
\]
It follows from Lemma~\ref{lem:lll} that
\[
\PP\left(
\bigcap_{v\in V(D)}\overline{E_v}
\cap
\bigcap_{\substack{v\in V(D)\\x_v>0}}\overline{G_v}
\right)>0.
\]

Hence there exists \(J\) such that none of the events \(E_v\) and \(G_v\) occurs.
For every \(v\in V(D)\setminus J\), the event \(E_v\) does not occur, so
\(N_D^+(v)\cap J\ne\varnothing\). Thus \(J\) is absorbing.

If \(v\in J\), then \(p_v>0\), and hence \(x_v>0\). Since \(G_v\) does not
occur, we have
\(d_{D[J]}^+(v)=X_v\le t-1\). Therefore,
\(\Delta^+(D[J])\le t-1\).

For every \(U\subseteq J\), the underlying graph of \(D[U]\) has at most
\((t-1)|U|\) edges and hence the average degree is at most \(2(t-1)\).
It follows that the underlying graph of \(D[J]\) is \(2(t-1)\)-degenerate.
Consequently, \(\chi(D[J])\le2t-1\).
\end{proof}

\vskip 2mm
Our second step is trying to remove the in-degree assumption from Lemma \ref{lem:bounded}. 

Let \(A=(A_{uv})_{u,v\in V(D)}\) be the adjacency matrix of \(D\), where \(A_{uv}=1\) if \(u\to v\), and \(A_{uv}=0\) otherwise. The proof uses a convergent matrix series associated with \(A\) to obtain positive weights \(w_v\) satisfying $\sum_{u\in N_D^-(v)}w_u<(d+1)w_v$ for every $v\in V(D)$. The
layers defined by $\ell(v)=\lfloor\log_{d+1}w_v\rfloor$ have bounded
in-degree, and an arc decreases the layer index by at most one. Within each
parity class, the layers can therefore be processed from high to low so that
the absorbing sets selected in distinct layers are nonadjacent. Their colours
can be reused, and two palettes suffice.

The next lemma partitions a bounded-out-degree digraph into layers of bounded
in-degree. Within each parity class, every arc between distinct layers is
directed from the lower layer to the higher one.

\begin{lemma}
\label{lem:layers}
Let \(d\ge1\) be an integer, and let \(D\) be an oriented graph with
\(\Delta^+(D)\le d\). There is a map \(\ell:V(D)\to\Zge\) such that, with
\(V_j=\{v:\ell(v)=j\}\),
\[
  u\to v\ \Longrightarrow\ \ell(v)\ge\ell(u)-1,
  \qquad
  \Delta^+(D[V_j])\le d,
  \quad
  \Delta^-(D[V_j])\le(d+1)^2-1
  \quad(j\ge0).
\]
In particular, if \(u\to v\), \(\ell(u)\equiv\ell(v)\pmod2\), and
\(\ell(u)\ne\ell(v)\), then \(\ell(u)<\ell(v)\).
\end{lemma}

\begin{proof}
Put $b=d+1$. We may assume that \(D\) is nonempty. For a matrix \(M\), let \(\|M\|_1\) denote its maximum absolute column sum. By the definition of adjacency matrix \(A\), the column sum of \(A^{\mathsf T}\) corresponding to a vertex \(v\) is \(d_D^+(v)\). Hence \(\|b^{-1}A^{\mathsf T}\|_1\le d/b<1\), and therefore the Neumann series
\[
w:=\sum_{r=0}^{\infty}b^{-r}(A^{\mathsf T})^r\one
\]
converges. Here \(\one\) is the all-ones vector in \(\mathbb{R}^{V(D)}\). Since its first term is \(\one\), we have \(w\ge\one\). Moreover, \(A^{\mathsf T}w=b(w-\one)<bw\).

Set \(\ell(v)=\lfloor\log_b w_v\rfloor\). Since
\[
  (A^{\mathsf T}w)_v=\sum_{u\in N_D^-(v)}w_u,
\]
an arc \(u\to v\) implies
\[
  b^{\ell(u)}\le w_u\le(A^{\mathsf T}w)_v
  <bw_v<b^{\ell(v)+2}.
\]
Thus \(\ell(u)\le\ell(v)+1\).

Fix \(v\in V_j\). Every in-neighbour of \(v\) in \(V_j\) has weight at
least \(b^j\), whereas \(w_v<b^{j+1}\). Hence
\[
  d^-_{D[V_j]}(v)b^j
  \le \sum_{\substack{u\to v\\u\in V_j}}w_u
  \le (A^{\mathsf T}w)_v
  <bw_v<b^{j+2}.
\]
Thus \(d^-_{D[V_j]}(v)<b^2\), and therefore
\(d^-_{D[V_j]}(v)\le b^2-1=(d+1)^2-1\). The out-degree bound follows from
\(\Delta^+(D)\le d\).

Finally, if \(u\to v\) and \(\ell(u)\equiv\ell(v)\pmod2\), then
\(\ell(v)-\ell(u)\) is even and at least \(-1\), and so is nonnegative. Hence \(\ell(v)>\ell(u)\) whenever the layers are distinct.
\end{proof}

The following lemma shows how to combine absorbing sets chosen from a sequence of layers in which every arc between distinct layers is directed from the lower layer to the higher one, without increasing the chromatic bound.

\begin{lemma}
\label{lem:assembly}
Let \(H\) be an oriented graph, let \(r,m\ge1\), and let \(V_1,\ldots,V_r\) be an ordered partition of
\(V(H)\) such that every arc between \(V_i\) and \(V_j\), where \(i<j\),
is directed from \(V_i\) to \(V_j\). Suppose that for every \(i\in[r]\) and
every \(U\subseteq V_i\), there exists an absorbing set \(S\subseteq U\) in
\(H[U]\) such that \(\chi(H[S])\le m\). Then \(H\) has an absorbing set
\(K\) such that \(\chi(H[K])\le m\).
\end{lemma}

\begin{proof}
By hypothesis, choose an absorbing set \(K_r\) of \(H[V_r]\) with
\(\chi(H[K_r])\le m\). Now let \(i<r\), and suppose that \(K_{i+1},\ldots,K_r\) have been chosen.
Set \(K_{>i}:=\bigcup_{j>i}K_j\) and
\(U_i:=\{v\in V_i:N_H^+(v)\cap K_{>i}=\varnothing\}\).
Thus \(U_i\) consists precisely of the vertices of \(V_i\) that are not
already absorbed by \(K_{>i}\). If \(U_i=\varnothing\), set
\(K_i=\varnothing\); otherwise, by hypothesis, choose an absorbing set
\(K_i\) of \(H[U_i]\) with \(\chi(H[K_i])\le m\).

Since every arc between \(V_i\) and \(V_j\) with \(j>i\) is directed from
\(V_i\) to \(V_j\), there is no arc from \(K_{>i}\) to \(U_i\), while the
definition of \(U_i\) ensures that there is no arc from \(U_i\) to
\(K_{>i}\). Hence there are no arcs between \(K_i\) and \(K_{>i}\).
It follows that there are no arcs between any two distinct sets among
\(K_1,\ldots,K_r\), and hence they may all be coloured using the same
palette of \(m\) colours.

Moreover, every vertex in \(V_i\setminus U_i\) has an out-neighbour in
\(K_{>i}\), while every vertex in \(U_i\setminus K_i\) has an out-neighbour
in \(K_i\). Hence \(K=\bigcup_{i=1}^rK_i\) is an absorbing set of \(H\),
and \(\chi(H[K])\le m\).
\end{proof}

\vskip 2mm
We are now in a position to to prove Theorem~\ref{thm:absorb}. For convenience, define \(a(0)=1\) and \(a(d)=4\left\lceil3\ln\bigl(2e(d+1)^3\bigr)\right\rceil-2\) for \(d\ge1\).

\begin{proof}[\bfseries{Proof of Theorem~\ref{thm:absorb}}]
We prove the stronger bound \(\chi(D[S])\le a(d)\). If \(d=0\), then \(D\) has no arcs. Taking \(S=V(D)\), we have that \(S\) is absorbing and \(\chi(D[S])\le 1=a(0)\). Assume that \(d\ge1\), and take the
layers from Lemma~\ref{lem:layers}.

Set \(h=(d+1)^2-1\),
\(t=\lceil3\ln\bigl(2e(d+1)^3\bigr)\rceil\), and \(m=2t-1\). For every
\(j\ge0\) and \(U\subseteq V_j\), the graph \(D[U]\) has maximum out-degree
at most \(d\) and maximum in-degree at most \(h\). Since
\((d+1)(h+1)=(d+1)^3\), by Lemma~\ref{lem:bounded}, we obtain an absorbing set of
chromatic number at most \(m\) in \(D[U]\).

For \(p\in\{0,1\}\), let
\(W_p:=\bigcup_{j\equiv p\pmod2}V_j\).
We construct an absorbing set \(K_p\) of \(D[W_p]\) with
\(\chi(D[K_p])\le m\). If \(W_p=\varnothing\), take \(K_p=\varnothing\).
Otherwise, consider the nonempty layers \(V_j\) with \(j\equiv p\pmod2\),
ordered by increasing \(j\). By Lemma~\ref{lem:layers}, if
\(j<j'\) and both \(V_j\) and \(V_{j'}\) are nonempty, then every arc
between \(V_j\) and \(V_{j'}\) is directed from \(V_j\) to \(V_{j'}\). Hence by Lemma~\ref{lem:assembly}, we obtain the
desired absorbing set \(K_p\).

Set \(K=K_0\cup K_1\). Every vertex outside \(K\) is absorbed within its own
parity class, and hence \(K\) is absorbing. Moreover,
\[
\chi(D[K])\le \chi(D[K_0])+\chi(D[K_1])
\le 2m=4t-2=a(d).
\]

For \(d\ge1\), using \(\lceil x\rceil<x+1\), we obtain
\[
a(d)<36\ln(d+1)+14+12\ln2.
\]
Define
\(
g(x):=31\log(x+2)-36\ln(x+1)-14-12\ln2.
\)
A direct calculation shows that \(g(5)>0\), while
\[
g'(x)=\frac{31}{(x+2)\ln2}-\frac{36}{x+1}>0
\qquad\text{for }x\ge5.
\]
Thus \(g(x)>0\) for all \(x\ge5\), and hence
\[
a(d)<31\log(d+2)
\]
for every \(d\ge5\). The cases \(1\le d\le4\) are verified directly from
\(
(a(1),a(2),a(3),a(4))=(46,58,70,78).
\)
\end{proof}

\section*{Acknowledgement}
\noindent This research is supported by the National Key R\&D Program of China under grant number 2024YFA1013900 and NSFC under grant number 12471327.

%\vspace{-1.2em}
\section*{Declaration}
\noindent\textbf{Conflict of interest.}
The authors declare that they have no known competing financial interests or personal relationships that could have appeared to influence the work reported in this paper.

%\vspace{-1.2em}
\section*{Data availability}
No data were used for the research described in this paper.

\end{document}